\documentclass{ifacconf}

\usepackage{graphicx}      
\usepackage{natbib}        
\usepackage{amsmath,amsxtra, amssymb,fancyhdr,latexsym}
\usepackage{bbm}
\begin{document}
\newcommand{\supp}{\operatorname*{supp}}
\newcommand{\pair}[2]{\left\langle #1 , #2 \right\rangle}
\newcommand{\Int}[4]{\int_{#1}^{#2}\! #3 \, #4}
\newcommand{\pd}[2]{\dfrac{\partial #1}{\partial #2}}
\newcommand{\map}[3]{#1 : #2 \rightarrow #3}
\def\R{\mathbb{R}}
\def\Rp{\mathbb{R}^+}
\def\N{\mathbb{N}}
\def\Np{\mathbb{N}^+}

\def\vd{v_{\text{d}}} 
\def\vs{v_{\text{s}}}
\def\vhat{\hat{v}}
\def\M{\mathcal{M}}
\def\tg{\tilde{g}}
\def\tgs{\tilde{g}_{\text{s}}}
\def\tga{\tilde{g}_{\text{a}}}
\def\sst{s_{\text{s}}}
\def\sat{s_{\text{a}}}
\newcommand{\g}[1]{g\left(\dfrac{#1}{|#1|}\cdot\dfrac{\vd}{|\vd|}  \right)}

\begin{frontmatter}

\title{Crowds reaching targets by maximizing entropy: a Clausius-Duhem inequality approach\thanksref{footnoteinfo}}

\thanks[footnoteinfo]{JE kindly acknowledges the financial support of the Netherlands Organisation for Scientific Research (NWO), Graduate Programme 2010.}

\author[First]{Joep Evers}
\author[Second]{Adrian Muntean}
\author[Third]{Fons van de Ven}

\address[First]{ICMS \& CASA, Eindhoven University of Technology,
   P.O. Box 513, 5600 MB Eindhoven, The Netherlands (e-mail: j.h.m.evers@tue.nl).}
\address[Second]{ICMS \& CASA, Eindhoven University of Technology,
   P.O. Box 513, 5600 MB Eindhoven, The Netherlands (e-mail: a.muntean@tue.nl)}
\address[Third]{CASA, Eindhoven University of Technology,
   P.O. Box 513, 5600 MB Eindhoven, The Netherlands (e-mail: a.a.f.v.d.ven@tue.nl)}

\begin{abstract}                
In this paper we propose the use of concepts from thermodynamics in the study of crowd dynamics. Our continuous model consists of the continuity equation for the density of the crowd and a kinetic equation for the velocity field. The latter includes a nonlocal term that models interactions between individuals. To support our modelling assumptions, we introduce an inequality that resembles the Second Law of Thermodynamics, containing an entropy-like functional. We show that its time derivative equals a positive dissipation term minus a corrector term. The latter term should be small for the time derivative of the entropy to be positive. In case of isotropic interactions the corrector term is absent. For the anisotropic case, we support the claim that the corrector term is small by simulations for the corresponding particle system. They reveal that this term is sufficiently small for the entropy still to increase. Moreover, we show that the entropy converges in time towards a limit value.
\end{abstract}

%

\end{frontmatter}
\section{Introduction}
Studying the behaviour of people in a crowd is nowadays no longer simply an activity of psychologists and social scientists. During the last few decades, it became clear that understanding and predicting the dynamics of people moving around is extremely important in the design of buildings and infrastructure, and for the safety of visitors at a large-scale event.\\
\\
Physicists and mathematicians started with providing answers to the questions that were posed. Their approach was (and is) similar to the way in which they deal with the physical, non-living world around us. Inevitably, this yields models in which people are treated as non-living material obeying physical laws. Illustrative for this way of thinking is the \textit{social force model}, see e.g.~\cite{Helbing95}. We remark that more or less in parallel and in the same spirit, the study of vehicular traffic developed. The aim of this paper is not to give an extensive overview of the existing models and literature. Instead, for more details the reader is referred e.g. to the excellent overview \cite{Bellomo} and the references cited therein.\\
\\
The difficulty in these models is in the constitutive relations. In physics, the intuition for these would be provided by experiments. We do not claim that doing experiments with pedestrians is impossible, but in any case it is difficult. One of the issues is the reproducibility.\\
\\
In this paper we explore a different way of justifying our modelling assumptions. We propose that a crowd obeys an inequality like the Second Law of Thermodynamics (also called Clausius-Duhem inequality) and thus maximizes `entropy'. If this is shown to be true, then this increases the consistency and trustworthiness of our model. For more background on thermodynamic concepts and their context, the reader is referred e.g.~to \cite{MullerRuggeri}.\\
\\
In Section \ref{sec: model} we describe our first-order continuum model including the constitutive relation for the velocity field. Next, we introduce in Section \ref{sec: inequality} a concept of (generalized) entropy. If all interactions are isotropic with respect to the direction in which another individual is perceived, then one can show that the time derivative of the entropy is non-negative. More work is required in the anisotropic case. Our numerical illustration in Section \ref{sec: numerics} suggests that even in that case the entropy inequality holds. This is the main conclusion here. We close the paper with an outlook on future work.

\section{Model equations}\label{sec: model}
Consider the continuity equation
\begin{equation}\label{continuity equation}
\pd{\rho}{t}+\nabla\cdot(\rho v)=0
\end{equation}
on $\R^d\times\Rp$ ($d\in\Np$ fixed). Our model is in the spirit of \cite{CristianiPiccoliTosin}, be it that their model is formulated in a more general way. We have a \textit{first-order} model, since we prescribe our velocity directly; see \cite{Coscia} for an exposition of first-order models \textit{versus} second-order models. We assume that the velocity field is the sum of two contributions
\begin{equation}\label{decomposition v}
v := \vd + \vs,
\end{equation}
a \textit{desired} velocity $\vd$, and a \textit{social} velocity $\vs$. The desired velocity is the velocity of an individual that is alone in the space $\R^d$. We take this velocity to be a constant, while
\begin{equation}\label{def vs}
\vs(x) := \Int{\R^d}{}{\g{x-y}\nabla W(|x-y|)\rho(y)}{dy}.
\end{equation}
The term $\vs$ models the interactions between individuals. Note that we often omit the explicit time dependence of our variables. Here, $\map{W}{\Rp}{\R}$ is the potential governing the interactions. We use the word `potential' here, since first-order models can be viewed as overdamped limits of second-order models. In the latter, $\nabla W$ corresponds to a (generalized) force and as such, $W$ is a potential. We refer to \cite{Edelstein2003}, p.~360, for a derivation of the overdamped limit in terms of a particle system.\\
In \eqref{def vs}
\begin{equation}\label{def nabla W}
\nabla W(|x-y|) := W'(|x-y|)\dfrac{x-y}{|x-y|}.
\end{equation}
The function $\map{g}{[-1,1]}{[0,1]}$ incorporates anisotropy in the model. This anisotropy arises because people have front and back sides. It depends on the direction in which one person perceives other people, how much influence they have on his motion. We restrict ourselves to linear functions $g$; that is, linear in $\frac{x-y}{|x-y|}\cdot\frac{\vd}{|\vd|}$, which is ($-1$ times) the cosine of the angle under which point $x$ sees point $y$. See \cite{Gulikers} for our previous investigations on the effect of $g$ on the dynamics of the underlying particle system.\\
\\
\cite{Coscia} do not focus on the nonlocal dependence on $\rho$. Unlike in \eqref{def vs}, they only allow for local dependence on $\rho$ and/or $\nabla\rho$.\\
\\
The choice of this velocity field is an ansatz. A way to justify this modelling assumption, is by proving that the model is somehow consistent with ideas from thermodynamics. This will be the main motivation for all steps in the sequel.\\
\\
We define the total mass of the crowd by
\begin{equation}
\M := \Int{\R^d}{}{\rho(x)}{dx},
\end{equation}
its centre of mass by
\begin{equation}
x_0 := \dfrac{1}{\M}\Int{\R^d}{}{x\rho(x)}{dx},
\end{equation}
and the velocity of its centre of mass (or: \textit{barycentric velocity}) by
\begin{equation}
v_0 := \dfrac{dx_0}{dt} = \dfrac{1}{\M}\Int{\R^d}{}{v(x)\rho(x)}{dx}.
\end{equation}
Finally, we introduce the velocity with respect to the velocity of the centre of mass as
\begin{equation}
\vhat(x) := v(x)-v_0.
\end{equation}

\section{Clausius-Duhem-like inequality}\label{sec: inequality}
Before introducing an entropy-like functional, and the corresponding Clausius-Duhem-like inequality, we define the dissipation function $\map{D}{[0,T]}{\Rp}$ (with $T>0$ some fixed final time)
\begin{equation}\label{def D}
D(t) := \Int{\R^d}{}{\rho(x)|\vhat(x)|^2}{dx}.
\end{equation}
This dissipation function can be derived from the second-order model (as in fluid mechanics) by multiplying the momentum equation with $v$. What we introduce here is a first-order version of this function, and therefore we keep the name \textit{dissipation function}. Justification of this choice is that variation of $D$ with respect to $\rho$ yields the kinetic equations \eqref{decomposition v}--\eqref{def vs}.\\
Since $\vd$ is constant,
\begin{align}
\nonumber v_0 =& \dfrac{1}{\M}\Int{\R^d}{}{v(z)\rho(z)}{dz}\\
\nonumber =& \vd\dfrac{1}{\M}\Int{\R^d}{}{\rho(z)}{dz} + \dfrac{1}{\M}\Int{\R^d}{}{\vs(z)\rho(z)}{dz}\\
=& \vd + \dfrac{1}{\M}\Int{\R^d}{}{\vs(z)\rho(z)}{dz},
\end{align}
so
\begin{equation}\label{characterization v-v0}
v(x)-v_0 = \vs(x)-\dfrac{1}{\M}\Int{\R^d}{}{\vs(z)\rho(z)}{dz}.
\end{equation}
The second term on the right-hand side of \eqref{characterization v-v0} is independent of $x$. For the dissipation, we then obtain
\begin{align}
\nonumber D(t) =& \Int{\R^d}{}{\rho\,(v-v_0)\cdot\vhat}{dx}\\
\nonumber =& \Int{\R^d}{}{\rho\,\vs\cdot\vhat}{dx} - \dfrac{1}{\M}\Int{\R^d}{}{\vs\rho}{dz}\cdot\Int{\R^d}{}{\rho\vhat}{dx}\\
=& \Int{\R^d}{}{\rho\,\vs\cdot\vhat}{dx},\label{D in terms of vsoc}
\end{align}
since $\Int{\R^d}{}{\rho\vhat}{dx}=0$, which follows from the definition of $\vhat$. For the ease of notation, we define
\begin{equation}
\tg(\xi):=\g{\xi},
\end{equation}
for all $\xi\in\R^d\setminus\{0\}$. Replacing $x$ by $y$ (and \textit{vice versa}) in $D$, we obtain
\begin{align}
\nonumber & D(t) =\\
\nonumber & \Int{\R^d}{}{\rho(x)\,\left(\Int{\R^d}{}{\tg(x-y)\nabla W(|x-y|)\rho(y)}{dy}\right)\cdot\vhat(x)}{dx}\\
\nonumber =& \Int{\R^d}{}{\rho(y)\,\left(\Int{\R^d}{}{\tg(y-x)\nabla W(|y-x|)\rho(x)}{dx}\right)\cdot\vhat(y)}{dy}\\
=& \Int{\R^d}{}{\rho(x)\Int{\R^d}{}{\tg(y-x)\nabla W(|y-x|)\cdot\vhat(y)\rho(y)}{dy}}{dx},
\end{align}
by changing the order of integration in the last step. We conclude from \eqref{def nabla W} that
\begin{equation}\label{nabla W x-y y-x}
\nabla W(|y-x|) = -\nabla W(|x-y|).
\end{equation}
Thus,
\begin{align}
\nonumber & D(t) =\\
& -\Int{\R^d}{}{\rho(x)\Int{\R^d}{}{\tg(y-x)\nabla W(|x-y|)\cdot\vhat(y)\rho(y)}{dy}}{dx}.\label{D g y-x}
\end{align}
A combination of \eqref{D in terms of vsoc} and \eqref{D g y-x} yields
\begin{equation}
2D(t) = \Int{\R^d}{}{\rho(x)\Int{\R^d}{}{\nabla W(|x-y|)\cdot V(x,y)\,\rho(y)}{dy}}{dx},
\end{equation}
where
\begin{align}
\nonumber V(x,y) :=& \vhat(x)\tg(x-y)-\vhat(y)\tg(y-x)\\
\nonumber =& (\vhat(x)-\vhat(y))\,\left[\dfrac12 \tg(x-y) + \dfrac12 \tg(y-x)\right]\\
\nonumber &+ (\vhat(x)+\vhat(y))\,\left[\dfrac12 \tg(x-y) - \dfrac12 \tg(y-x)\right]\\
=:& (\vhat(x)-\vhat(y))\tgs(x-y) + (\vhat(x)+\vhat(y))\tga(x-y).
\end{align}
Here, $\tgs$ and $\tga$ are the symmetric and antisymmetric parts of $\tg$, respectively. Thus
\begin{align}\label{D two terms v hat tilde g}
\nonumber D(t) =& \,\dfrac12\int_{\R^d}\!\rho(x)\int_{\R^d}\!\nabla W(|x-y|)\cdot\\
\nonumber &\hspace{2 cm}\cdot(\vhat(x)-\vhat(y))\tgs(x-y)\rho(y)\,dy\,dx\\
\nonumber &+\,\dfrac12\int_{\R^d}\!\rho(x)\int_{\R^d}\!\nabla W(|x-y|)\cdot\\
\nonumber &\hspace{2 cm}\cdot(\vhat(x)+\vhat(y))\tga(x-y)\rho(y)\,dy\,dx\\
=:&\sst(t)+\sat(t).
\end{align}
Inspired by \cite{CarrilloMoll}, we define the following entropy-like functional
\begin{equation}\label{def S}
S(t) := \dfrac12 \Int{\R^d}{}{\rho(x)\Int{\R^d}{}{W(|x-y|)\tg(x-y)\rho(y)}{dy}}{dx}.
\end{equation}
Since we have assumed that the function $g$ is linear, $\tgs$ is a constant, which we call $\alpha$. This brings us to
\begin{lem}
For linear $g$
\begin{align}
\nonumber&\dfrac{dS}{dt}=\\
&\dfrac\alpha 2\Int{\R^d}{}{\rho(x)\Int{\R^d}{}{\nabla W(|x-y|)\cdot(\vhat(x)-\vhat(y))\rho(y)}{dy}}{dx}.\label{dSdt lemma}
\end{align}
\end{lem}
\begin{pf}
\begin{align}
\nonumber \dfrac{dS}{dt}=&\dfrac12 \Int{\R^d}{}{\pd{\rho(x)}{t}\Int{\R^d}{}{W(|x-y|)\tg(x-y)\rho(y)}{dy}}{dx}\\
&+\dfrac12 \Int{\R^d}{}{\rho(x)\Int{\R^d}{}{W(|x-y|)\tg(x-y)\pd{\rho(y)}{t}}{dy}}{dx}.
\end{align}
Interchanging the order of integration in the first term, replacing $x$ by $y$ (and \textit{vice versa}) in the second term, and using $W(|y-x|)=W(|x-y|)$, we obtain
\begin{align}
\nonumber \dfrac{dS}{dt}=&\dfrac12 \Int{\R^d}{}{\rho(y)\Int{\R^d}{}{W(|x-y|)\tg(x-y)\pd{\rho(x)}{t}}{dx}}{dy}\\
\nonumber &+\dfrac12 \Int{\R^d}{}{\rho(y)\Int{\R^d}{}{W(|x-y|)\tg(y-x)\pd{\rho(x)}{t}}{dx}}{dy}\\
\nonumber =& \Int{\R^d}{}{\rho(y)\Int{\R^d}{}{W(|x-y|)\tgs(x-y)\pd{\rho(x)}{t}}{dx}}{dy}\\
=& -\alpha\Int{\R^d}{}{\rho(y)\Int{\R^d}{}{W(|x-y|)\nabla\cdot(\rho(x)v(x))}{dx}}{dy},
\end{align}
where, in the last step, we used the continuity equation \eqref{continuity equation}. Next, we apply integration by parts in $x$ (where we assume vanishing boundary terms) to obtain
\begin{align}
\nonumber \dfrac{dS}{dt}=& \alpha\Int{\R^d}{}{\rho(y)\Int{\R^d}{}{\nabla W(|x-y|)\cdot v(x)\rho(x)}{dx}}{dy}\\
\nonumber =& \dfrac\alpha 2\Int{\R^d}{}{\rho(y)\Int{\R^d}{}{\nabla W(|x-y|)\cdot v(x)\rho(x)}{dx}}{dy}\\
\nonumber & + \dfrac\alpha 2\Int{\R^d}{}{\rho(x)\Int{\R^d}{}{\nabla W(|y-x|)\cdot v(y)\rho(y)}{dy}}{dx}\\
\nonumber =& \dfrac\alpha 2\Int{\R^d}{}{\rho(x)\Int{\R^d}{}{\nabla W(|x-y|)\cdot v(x)\rho(y)}{dy}}{dx}\\
& - \dfrac\alpha 2\Int{\R^d}{}{\rho(x)\Int{\R^d}{}{\nabla W(|x-y|)\cdot v(y)\rho(y)}{dy}}{dx}.
\end{align}
In the last step we interchanged the order of integration in the first term, and used \eqref{nabla W x-y y-x} in the second term.\\
All this yields
\begin{align}
\nonumber &\dfrac{dS}{dt}=\\
& \dfrac\alpha 2\Int{\R^d}{}{\rho(x)\Int{\R^d}{}{\nabla W(|x-y|)\cdot (v(x)-v(y))\rho(y)}{dy}}{dx},
\end{align}
and the desired result then easily follows from the observation that
\begin{equation}
v(x)-v(y)= (\vhat(x)+v_0)-(\vhat(y)+v_0)=\vhat(x)-\vhat(y).
\end{equation}
$\square$\end{pf}
\textbf{Remark: }One can relax the linearity condition on $g$ and obtain the same result. It suffices to have $g'(\eta)=g'(-\eta)$ for all $\eta\in[-1,1]$. Then $\nabla\tgs\equiv 0$ and the corresponding term after integration by parts vanishes. For general (differentiable) $g$, this term remains, and thus an extra term appears in \eqref{dSdt lemma} (and eventually in the entropy inequality we are deriving). These technicalities are however beyond the scope and aim of this paper.\\
\\
The following equation summarizes our ideas so far in a condensed form:
\begin{equation}
\dfrac{dS(t)}{dt} =  D(t) - \sat(t),
\end{equation}
which is the more desirable form, as it leads us to an entropy-like inequality. We first observe (by its definition \eqref{def D}) that $D(t)\geqslant0$ for all $t$. Moreover, we identify the special case of fully isotropic interactions with $g\equiv1$. For $g\equiv1$ in \eqref{def vs}, the interaction term in $\vs$ is $\nabla W(|x-y|)$, which is the gradient of a radially symmetric function. In this case $\tgs\equiv\tg$ and $\tga\equiv0$. Thus, $\sat(t)=0$ for all $t$, and consequently
\begin{equation}
\dfrac{dS(t)}{dt} =  D(t) \geqslant 0,
\end{equation}
which is a Clausius-Duhem-type inequality. This is a special case of what was treated by \cite{CarrilloMoll}, although there the reduction by subtracting $v_0$ is not done. If $W$ is bounded, $S$ has an upper bound that is uniform in time. This implies that $S$ will tend to some limiting value as $t\rightarrow\infty$.\\
\\
Our main question is now whether similar conclusions can be drawn if $g$ is not constant. This is the case in which anisotropy is present in the interactions. In other words, some directions have more influence than others. In general, $\sat$ will then no longer be $0$. However, if this term is small (compared to $D$), still
\begin{equation}
\dfrac{dS(t)}{dt} \geqslant 0
\end{equation}
holds, which is the inequality we are looking for. In the sequel, we test this hypothesis numerically for a specific particle system.

\section{Numerical illustration of the anisotropic case}\label{sec: numerics}
\begin{figure}
\begin{center}
\includegraphics[width=8.4cm]{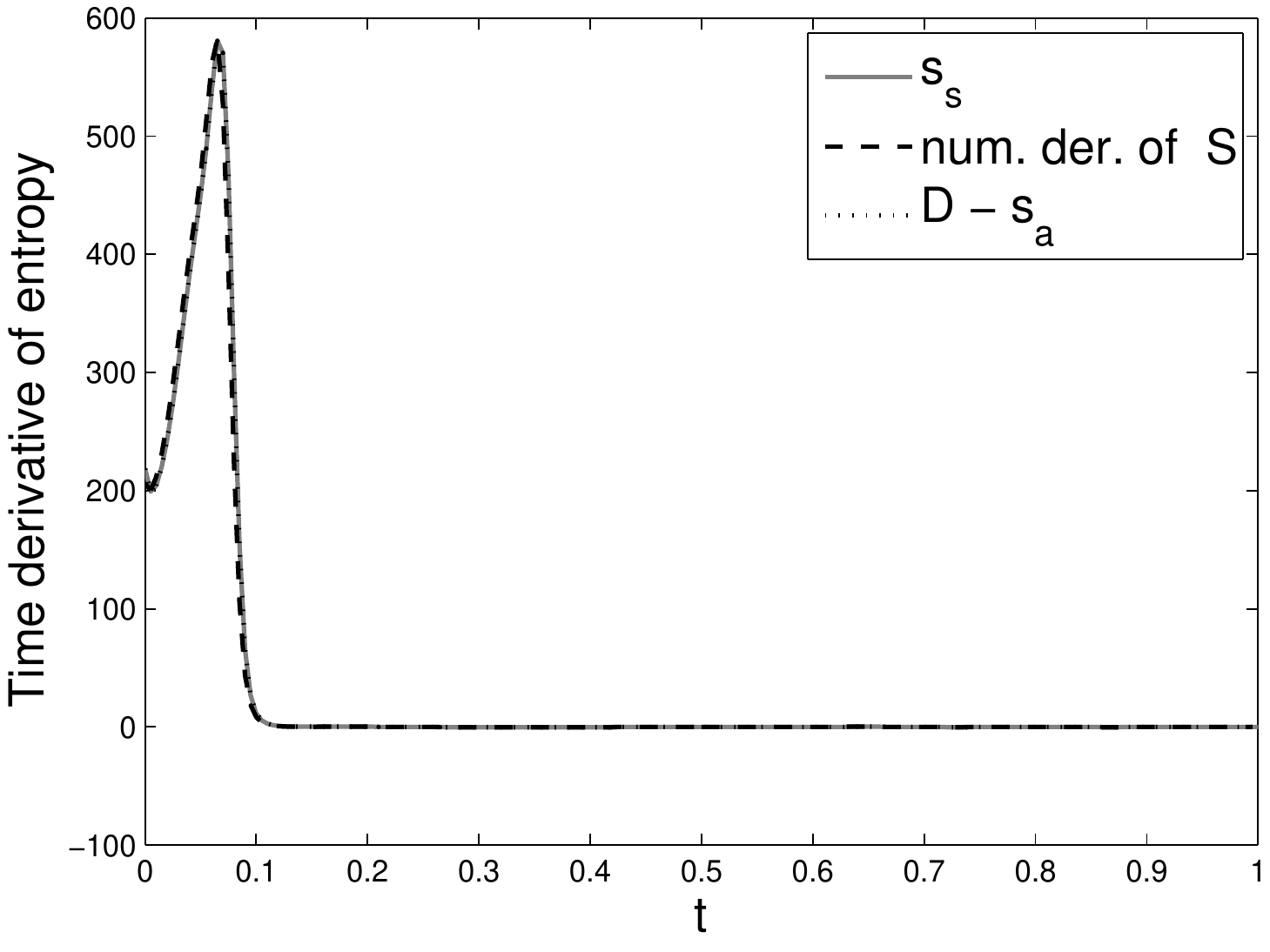}    
\caption{Plots of $\sst$, $dS/dt$ (numerical derivative of $S$), and $D-\sat$ as functions of time. Theory predicts that they should be the same, as is supported by the graphs. It is important to note that the curves are positive and decay to zero (fluctuations around zero are $\mathcal{O}(0.1)$). Results for a single simulation run.}
\label{fig:timeDerivative1}
\end{center}
\end{figure}

For a numerical illustration and investigation of the ideas described above, we simulate the particle system (of size $N$) corresponding to the model in Section \ref{sec: model}. In particular, we take $d=2$, and choose $W$ to be the \textit{Morse potential}
\begin{equation}
W(s) := C_a\,e^{-s/l_a}-C_r\,e^{-s/l_r},
\end{equation}
see \cite{Edelstein2003}, p.~363, or \cite{DOrsogna} and the references cited therein. Note that \cite{DOrsogna} use this potential in a second-order model. We demand that the parameters are positive and obey $l_r < l_a$ and $C_r/l_r > C_a/l_a$. This makes sure that the interactions are repulsive in the short range, and attractive in the long range. (NB: These conditions include the Case 4 mentioned in \cite{Edelstein2003}: $C_r>C_A$ and $l_a>l_r$.)\\
\\
Furthermore, we take
\begin{equation}
g(\eta) := \dfrac12(1+\sigma)-\dfrac12(1-\sigma)\eta.
\end{equation}
Here $\sigma\in[0,1]$ is a parameter that is used to tune the amount of anisotropy. The isotropic case corresponds to $\sigma=1$. Note that $\sigma$ relates to the aforementioned $\alpha$ via $\alpha=(1+\sigma)/2$. For our test, we take $\sigma=0.5$ and $N=25$. The former choice ensures us that we test the \textit{anisotropic} case. The latter choice reflects the current status of our work. Simulating the continuum is work in progress, as well as increasing the number of particles, as is mentioned in the outlook in Section \ref{sec: conclusions} (cf.~also \cite{Gulikers} in this respect).\\
\\
Initially the particles (individuals) are distributed randomly over the unit square $[0,1]^2$. We first show results for one simulation run, and secondly for a sequence of $1000$ runs. New initial conditions are generated in each run.\\
\\
To see whether our numerics comply with $dS/dt\geqslant0$, we show in Fig.~\ref{fig:timeDerivative1} the evolution in time of $dS/dt$. We deduced this quantity in three different ways from our simulations: from $\sst$, from the time derivative of $S$ (where $S$ was calculated for the particle system and the numerical time derivative was computed afterwards), and from $D-\sat$. Theoretically these are identical, and the plots show that this is also the case numerically. We will thus not bother about this issue any more in the sequel. The main conclusion from Fig.~\ref{fig:timeDerivative1} is that $dS/dt$ is indeed positive and, moreover, that it decays to zero.

\begin{figure}
\begin{center}
\includegraphics[width=8.4cm]{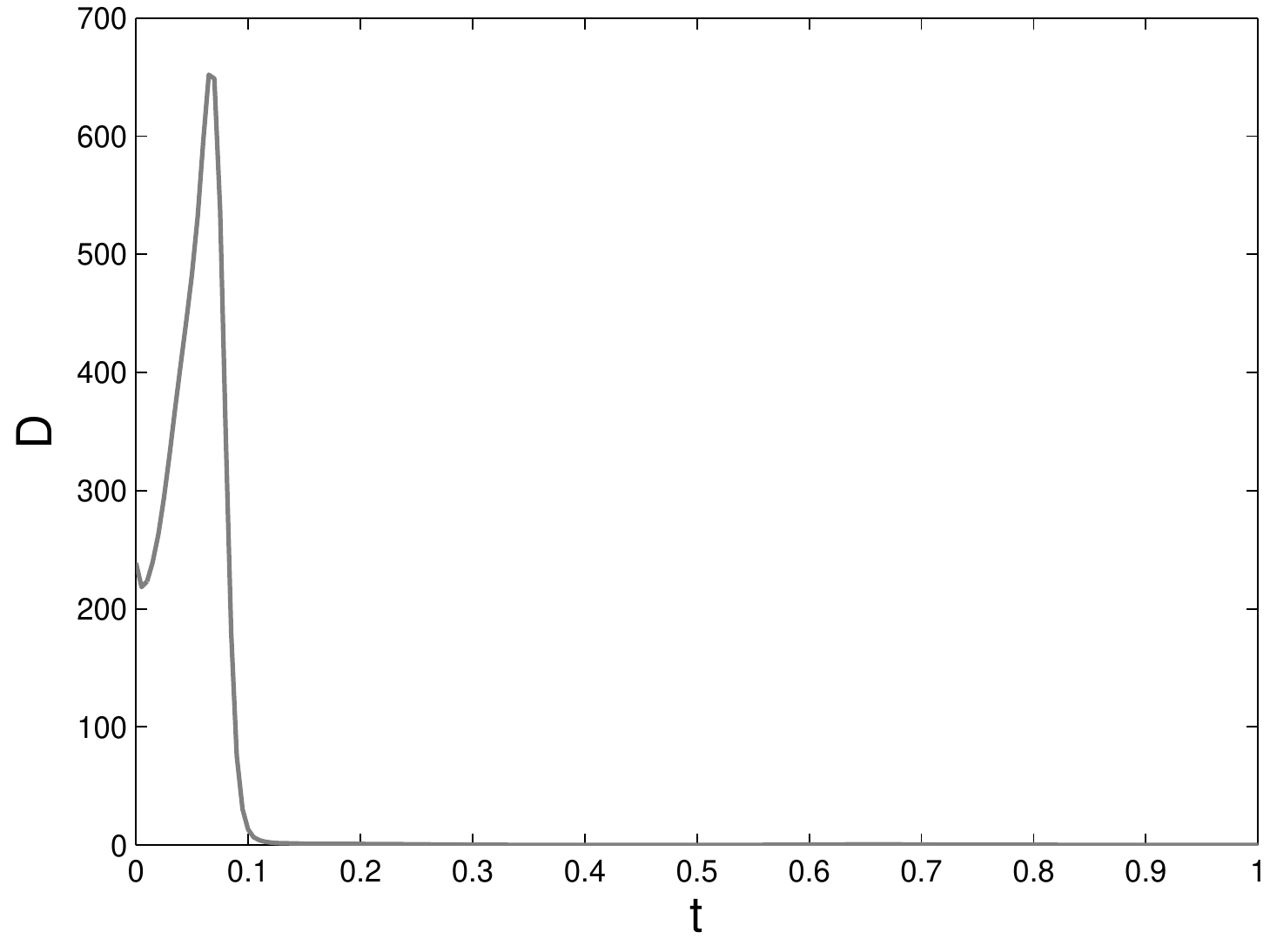}    
\caption{Plot of $D$ as a function of time. The value decays to zero. Results for a single simulation run.}
\label{fig:dissipation1}
\end{center}
\end{figure}

\noindent Next, we examine the behaviour in time of the dissipation $D$. By definition, $D$ is positive -- see \eqref{def D} -- but Fig.~\ref{fig:dissipation1} shows that $D$ tends to zero as $t$ increases. Note that the behaviour of $D$ is very similar to that of $dS/dt$.

\begin{figure}
\begin{center}
\includegraphics[width=8.4cm]{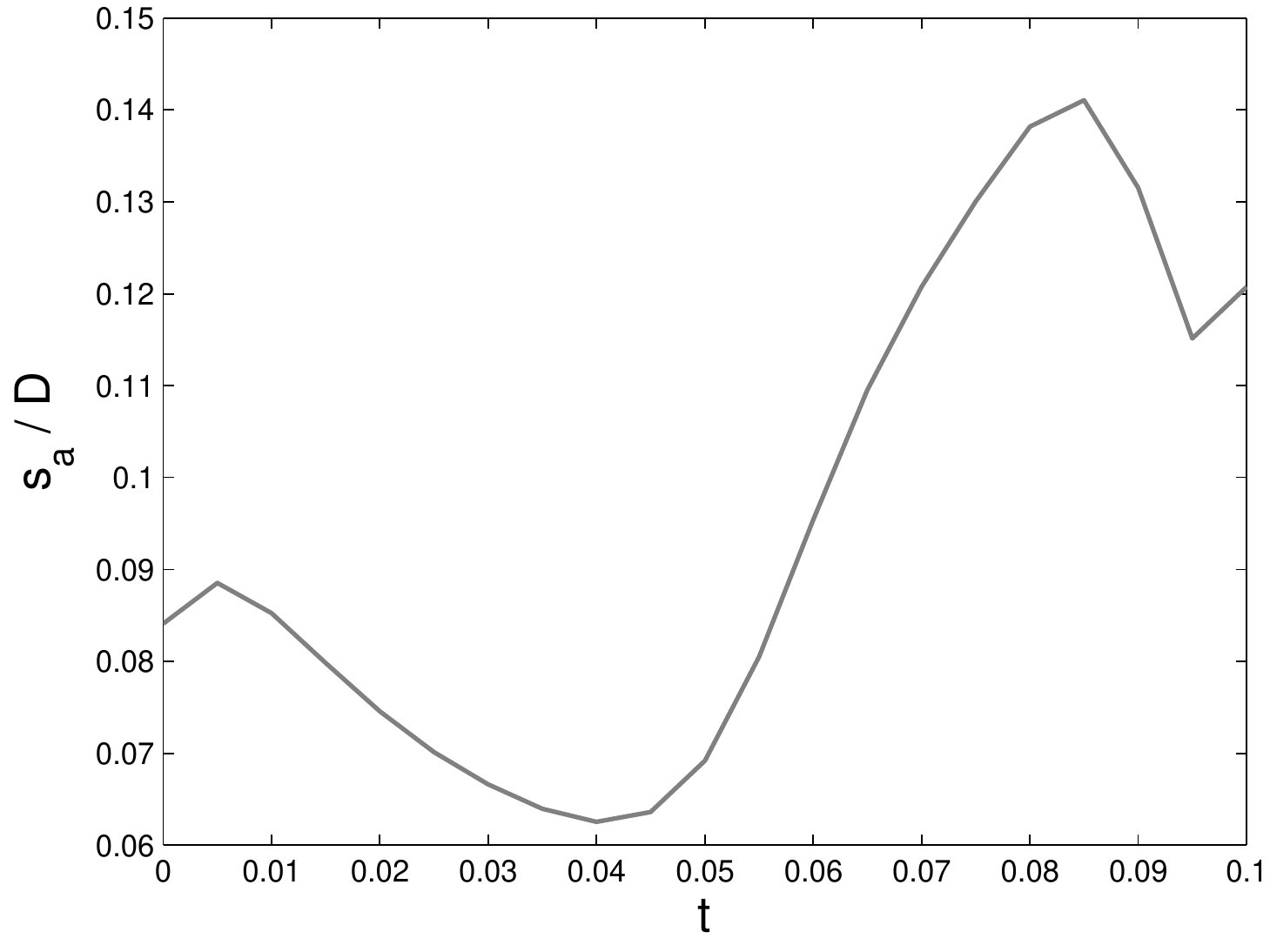}    
\caption{Plot of $\sat/D$ as a function of time on the interval $[0,0.1]$. The values are small. Results for a single simulation run.}
\label{fig:smallTerm1}
\end{center}
\end{figure}

\noindent A sharp transition is observed in Figs.~\ref{fig:timeDerivative1} and \ref{fig:dissipation1} for $t\approx0.1$ where the graphs become zero. If $D-\sat$ and $D$ are (nearly) zero, then $\sat$ must also be (nearly) zero. The hypothesis that $\sat$ is small, is thus valid after $t=0.1$. In Fig.~\ref{fig:smallTerm1}, we plotted the ratio $\sat/D$ on the interval $[0,0.1]$. We take this ratio since we are particularly interested in the size of $\sat$ \textit{with respect to $D$}. Indeed the values are small, i.e. of order $\mathcal{O}(0.1)$. We do not continue the graph after $t=0.1$ for a simple reason. Since both $\sat$ and $D$ are very small then, we would divide two small numbers to get the ratio $\sat/D$. The corresponding outcome does not provide any useful information.

\begin{figure}
\begin{center}
\includegraphics[width=8.4cm]{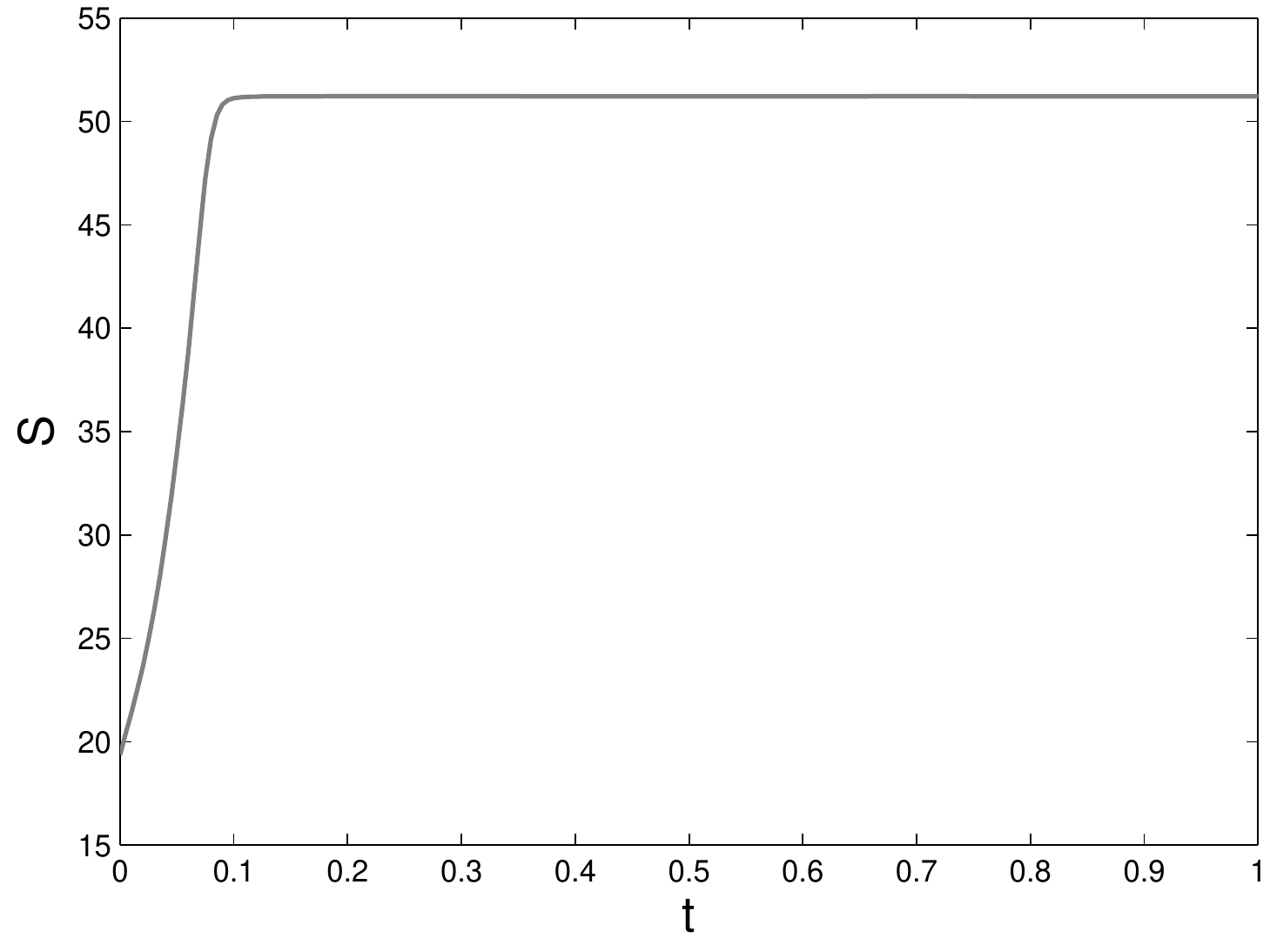}    
\caption{Plot of $S$ as a function of time. The value approaches a limit value. Results for a single simulation run.}
\label{fig:entropy1}
\end{center}
\end{figure}

\noindent Since the interaction potential $W$ is bounded, the entropy $S$ has a finite upper bound. As $dS/dt$ is positive (cf.~Fig.~\ref{fig:timeDerivative1}), we expect $S$ to approach a limit value, which is confirmed by Fig.~\ref{fig:entropy1}.

\begin{figure}
\begin{center}
\includegraphics[width=8.4cm]{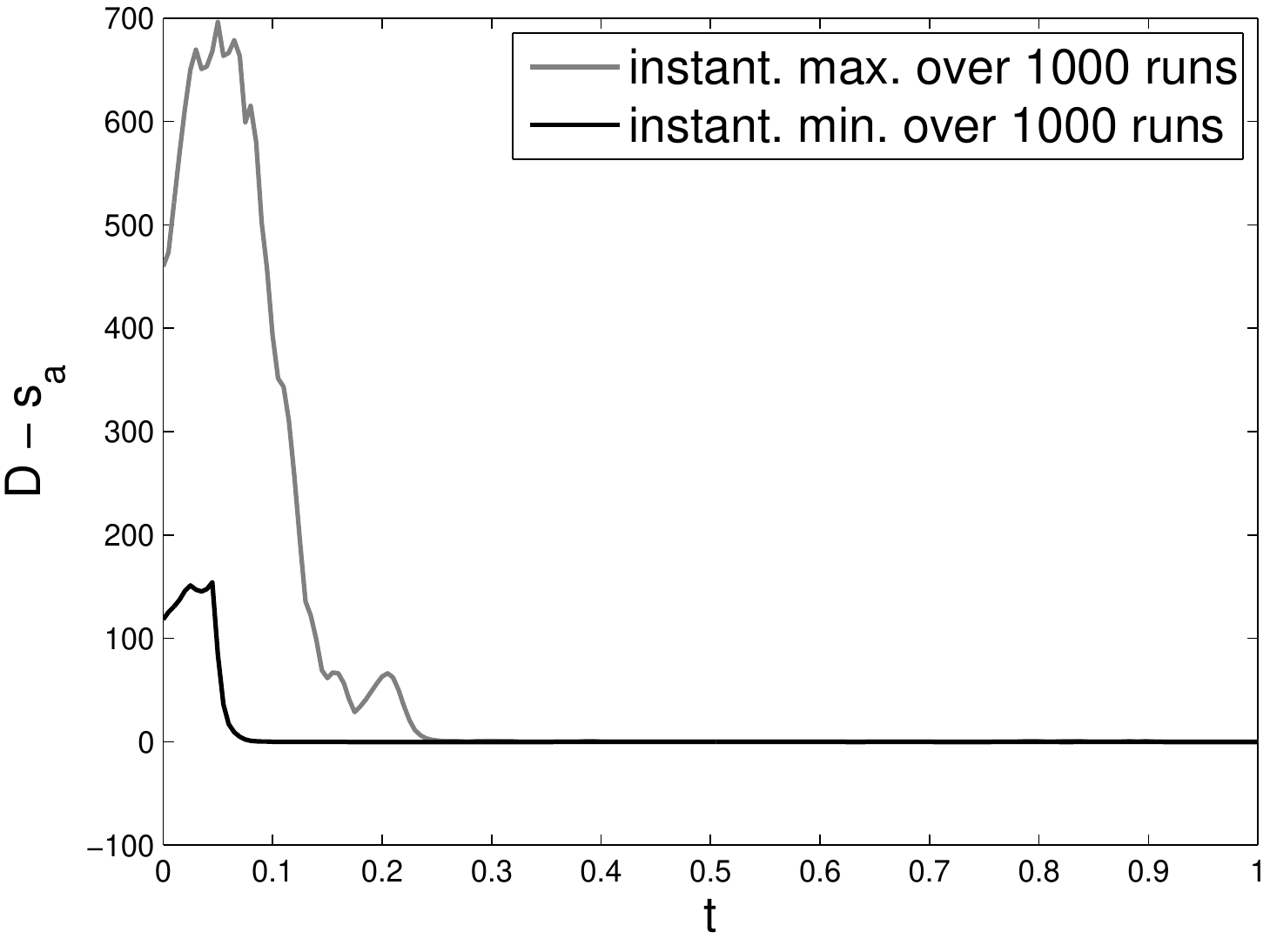}    
\caption{Plot of $D-\sat$ as a function of time. At every time instance both the maximum and the minimum are taken over 1000 simulation runs. The graphs support the claim that $D-\sat$ is positive and decays to zero.}
\label{fig:timeDerivative1000}
\end{center}
\end{figure}

\noindent We check now whether the conclusions from one simulation run also hold for multiple runs (where in each run we impose different initial conditions), to be sure that we have not just been `lucky' so far. In Fig.~\ref{fig:timeDerivative1000} we show two curves corresponding to $D-\sat$. At each time instance, we show both the minimum and maximum over all simulation runs. These graphs provide further evidence that $D-\sat$ (and hence $dS/dt$) is positive -- which follows from the minimum -- and decays to zero. Moreover, the term $\sat$ is thus small compared to $D$. After $t\approx0.25$ the deviation from zero of the two curves is of order $\mathcal{O}(0.1)$.

\begin{figure}
\begin{center}
\includegraphics[width=8.4cm]{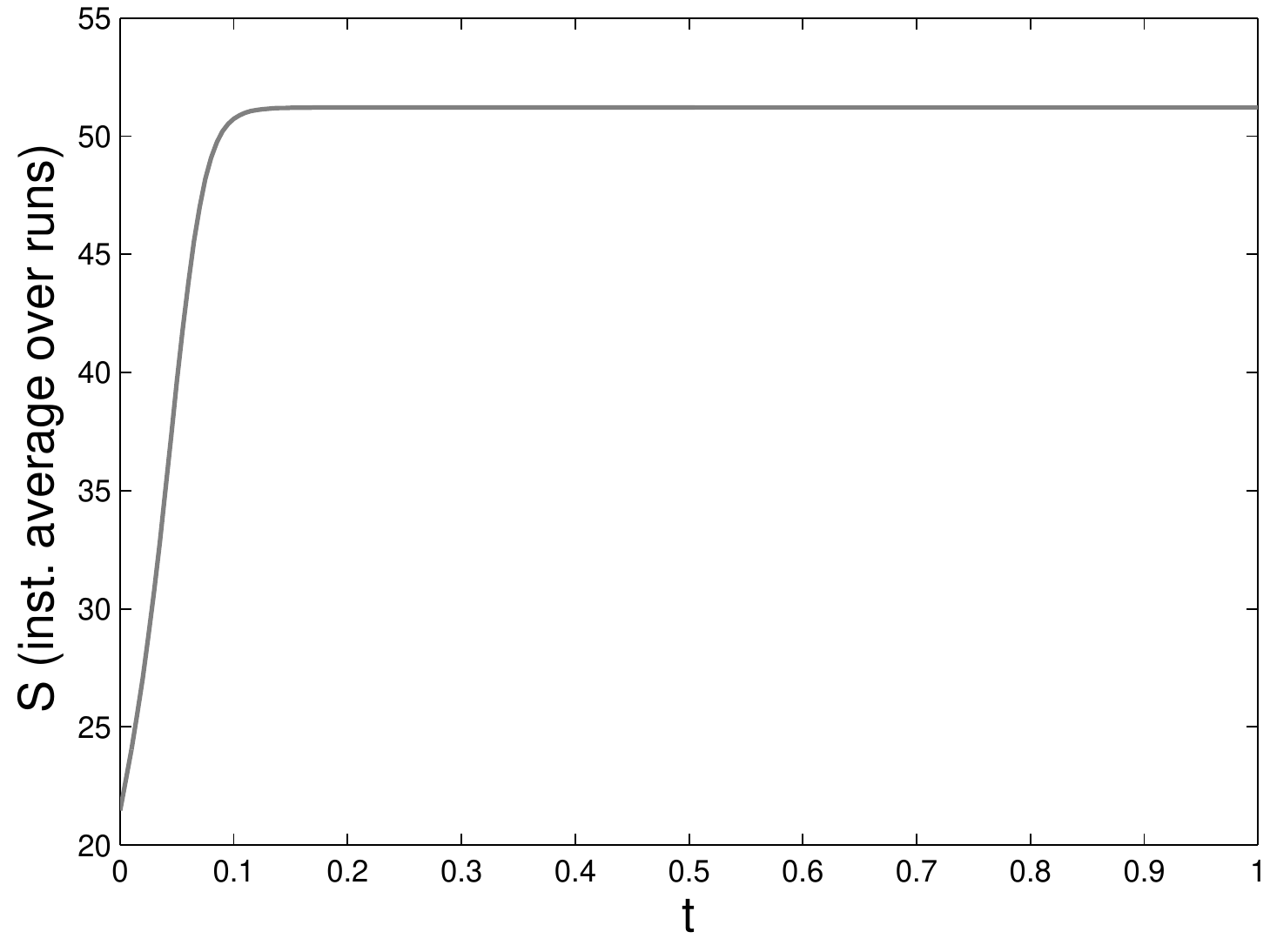}    
\caption{Plot of $S$ as a function of time. The value approaches a limit value. At every time instance the average is taken over 1000 simulation runs. We checked that in all runs the same limit value was attained.}
\label{fig:entropy1000}
\end{center}
\end{figure}

%
\noindent Taking into account Fig.~\ref{fig:timeDerivative1000}, which shows that $dS/dt\geqslant0$, we expect also $S$ to increase and eventually converge to a limit. In Fig.~\ref{fig:entropy1000} the average of $S$ over the $1000$ runs is shown, and this average indeed increases towards a limit. A plot (this graph is omitted here) of $\max_k |S_1(t)-S_k(t)|$ against time (where the index $k$ runs over all simulation instances, so $S_1$ is $S$ obtained from the first simulation run), decays to zero, implying that the limit value is the same in each simulation run.

\section{Conclusions and outlook}\label{sec: conclusions}
What does it mean for our system that the dissipation $D$ goes to zero? By its definition in \eqref{def D}, $D=0$ automatically implies that $\vhat\equiv0$. Since, for large $t$, the dissipation $D$ goes to zero, the system reaches a limit state, in which all material points have the same velocity (namely the velocity of the centre of mass $v_0(t)$). The consequence for the continuum model is that the density profile $\rho$ does not change shape any more. However, this does \textit{not} mean that the density $\rho$ is uniform. One can show that $\rho$ is conserved along characteristics $x(t)$ defined by the equation $dx(t)/dt=v_0(t)$. The configuration of the system is thus just convected.\\
When $D\equiv 0$, $S$ is constant in time. This follows, since $\vhat\equiv0$ implies $\sat\equiv0$, cf.~\eqref{D two terms v hat tilde g}, and thus $dS/dt\equiv0$.\\
The above explanation relates our entropy to the theory of Lyapunov functionals.\\
\\
The numerics presented in this paper suggest that even in the case of anisotropic interactions, the dynamics still obey a Clausius-Duhem-type inequality. This supports the idea that it is worthwhile to do more effort to prove this analytically. We are aware of the fact that this might only be true under certain technical conditions, which were satisfied (more or less `by accident'?) in our numerics. One of the aims of further theoretical investigations is to identify these conditions, and make them as sharp as possible. In contrast to the isotropic case, it is clear from our (unsuccessful) attempts up to now that it will not be an easy task to prove the inequality in the anisotropic case.\\
We stress here that proving an entropy inequality is not a goal in itself. It provides support for the `thermodynamic consistency' of our model, and as such tells us that our ansatz for the velocity field is admissible, and hence, not completely wrong.\\
\\
To get a better understanding of how to proceed towards these proofs, we suggest to do first some additional numerical experiments. In particular, we want to:
\begin{itemize}
  \item Do the same simulations, but change the value of $\sigma$. Note that $\tga(x-y)$ is proportional to $(1-\sigma)$, and this factor thus appears in $\sat$. The amplitude of this term most likely increases automatically with decreasing $\sigma$ ($\sigma\downarrow0$).\footnote{One should however be careful here in drawing this conclusion. Taking $\sigma$ closer to $0$ also changes the dynamics which, in turn, might cause a change in the integral term over $\rho$ in $\sat$. This change could counterbalance the change in $1-\sigma$ in such a way that our claim on the size of $\sat$ does not hold.} Following this line of argument, we expect that $dS/dt$ is only positive for $\sigma$ within a certain distance from $1$. This can lead to a condition on $\sigma$; possibly a condition in which $\sigma$ is combined with other quantities.
  \item Increase the number of particles $N$. One might guess that in a certain scaling, one can obtain the continuum model treated in this paper from the corresponding particle system that was used in the simulation section. We hope that for larger $N$ the outcome of the particle system will be closer to the result of the continuum model.
      In this paper, a relatively small $N$ was used, to get a system of ODEs that can still be handled easily. After all, the numerics in this paper are only intended to illustrate our ideas and confirm our conjectures. Considering the continuum limit and the limit process for $N\rightarrow\infty$ are subject of ongoing work. 
  \item To see whether our results depend on the precise choice of interaction potential $W$. The simplest way for this is to consider different parameters ($C_r, C_a, l_r, l_a$). However, it is much more interesting to try another type of potential (possibly with a singularity around the origin). If the repulsive behaviour around the origin is e.g.~of the type $\sim 1/r$ (cf.~Coulomb interactions), we lose the trivial upper bound on $S$ that followed from $\|W\|_{\infty}$. Numerics should then provide insight in whether the entropy still increases towards a limit value. Also, one could especially look at potentials that only model a zone of repulsion (that is, no attraction zone). This is interesting, since for repulsive interactions mass will completely `diffuse', and for each $x$ fixed, $\rho(x,t)\rightarrow0$ as $t\rightarrow\infty$. This implies that no steady states are to be expected. To what extent this destroys our entropy inequality is yet to be investigated.
\end{itemize}


\bibliography{references}             

\begin{thebibliography}{9}
\providecommand{\natexlab}[1]{#1}
\providecommand{\url}[1]{\texttt{#1}}
\providecommand{\urlprefix}{URL }
\expandafter\ifx\csname urlstyle\endcsname\relax
  \providecommand{\doi}[1]{doi:\discretionary{}{}{}#1}\else
  \providecommand{\doi}{doi:\discretionary{}{}{}\begingroup
  \urlstyle{rm}\Url}\fi

\bibitem[{Bellomo and Dogbe(2011)}]{Bellomo}
Bellomo, N. and Dogbe, C. (2011).
\newblock On the modelling of traffic and crowds: A survey of models,
  speculations, and perspectives.
\newblock \emph{SIAM Review}, 53(3), 409--463.

\bibitem[{Carrillo and Moll(2009)}]{CarrilloMoll}
Carrillo, J. and Moll, J. (2009).
\newblock Numerical simulation of diffusive and aggregation phenomena in
  nonlinear continuity equations by evolving diffeomorphisms.
\newblock \emph{SIAM J. Sci. Comput.}, 31, 4305--4329.

\bibitem[{Coscia and Canavesio(2008)}]{Coscia}
Coscia, V. and Canavesio, C. (2008).
\newblock First-order macroscopic modelling of human crowd dynamics.
\newblock \emph{Math. Mod. Meth. Appl. Sci.}, 18(suppl.), 1217--1247.

\bibitem[{Cristiani et~al.(2011)Cristiani, Piccoli, and
  Tosin}]{CristianiPiccoliTosin}
Cristiani, E., Piccoli, B., and Tosin, A. (2011).
\newblock Multiscale modeling of granular flows with application to crowd
  dynamics.
\newblock \emph{Multiscale Model. Simul.}, 9(1), 155--182.

\bibitem[{D'Orsogna et~al.(2006)D'Orsogna, Chuang, Bertozzi, and
  Chayes}]{DOrsogna}
D'Orsogna, M., Chuang, Y., Bertozzi, A., and Chayes, L. (2006).
\newblock Self-propelled particles with soft-core interactions: Patterns,
  stability, and collapse.
\newblock \emph{Phys. Rev. Lett.}, 96, 104302.

\bibitem[{Gulikers et~al.(2013)Gulikers, Evers, Muntean, and Lyulin}]{Gulikers}
Gulikers, L., Evers, J., Muntean, A., and Lyulin, A. (2013).
\newblock The effect of perception anisotropy on particle systems describing
  pedestrian flows in corridors.
\newblock \emph{J. Stat. Mech.}, P04025.

\bibitem[{Helbing and Moln\'ar(1995)}]{Helbing95}
Helbing, D. and Moln\'ar, P. (1995).
\newblock Social force model for pedestrian dynamics.
\newblock \emph{Phys. Rev. E}, 51(5), 4282--4286.

\bibitem[{Mogilner et~al.(2003)Mogilner, Edelstein-Keshet, Bent, and
  Spiros}]{Edelstein2003}
Mogilner, A., Edelstein-Keshet, L., Bent, L., and Spiros, A. (2003).
\newblock Mutual interactions, potentials, and individual distance in a social
  aggregation.
\newblock \emph{J. Math. Biol.}, 47, 353--389.

\bibitem[{M\"{u}ller and Ruggeri(1998)}]{MullerRuggeri}
M\"{u}ller, I. and Ruggeri, T. (1998).
\newblock \emph{{R}ational {E}xtended {T}hermodynamics}.
\newblock Springer Verlag.

\end{thebibliography}

\end{document}